\documentclass[preprint,12pt]{elsarticle}

\usepackage{graphicx}
\usepackage{amsmath}
\usepackage{amssymb}
\usepackage{hyperref}
\usepackage{xcolor}
\usepackage{mathtools}
\usepackage{subcaption}
\usepackage{booktabs}
\usepackage{bm}

\hypersetup{
    colorlinks=true,
    linkcolor=blue,
    filecolor=magenta,      
    citecolor=red,
    urlcolor=cyan,
}

\begin{document}

\begin{frontmatter}
\title{Sequential Piecewise Linear Programming for Convergent Optimization of Non-Convex Problems}
\author{James P.L. Tan}
\ead{jamestan@ntu.edu.sg}
\address{Nanyang Technological University, Singapore}
\begin{abstract}
A sequential piecewise linear programming method is presented where bounded domains of non-convex functions are successively contracted about the solution of a piecewise linear program at each iteration of the algorithm. Although feasibility and optimality are not guaranteed, we show that the method is capable of obtaining convergent and optimal solutions on a number of Nonlinear Programming (NLP) and Mixed Integer Nonlinear Programming (MINLP) problems using only a small number of breakpoints and integer variables. 
\end{abstract}
\end{frontmatter}

\section{Introduction}
Piecewise linear approximations are commonly used methods to approximate non-convex problems as mixed-integer programs \cite{Belotti1}. Such an approach has the benefit of being able to approximate a solution for non-convex problems while also being derivative-free and capable of global optimization \cite{Geissler1}. The main drawbacks being the introduction of additional integer variables that can considerably increase the computational complexity of the optimization problem, especially if the non-linear function is non-separable (i.e. cannot be written as a sum of one-dimensional functions). Also, if one wishes to get a better approximation of the solution, then even more breakpoints and integer variables have to be introduced, further increasing computational complexity. In this paper, we present SPPA (Sequential Piecewise Planar Approximation), a simple scheme of iteratively contracting the bounded domains about the solutions of a piecewise linear program. We show that with a low number of breakpoints, it is sufficient to converge on the optimal solution of a majority of MINLP problems we have tested. 

\subsection{Decomposing nonlinear functions into piecewise linear segments}
\begin{figure}[h!]
\centering
\includegraphics[width=0.8\textwidth]{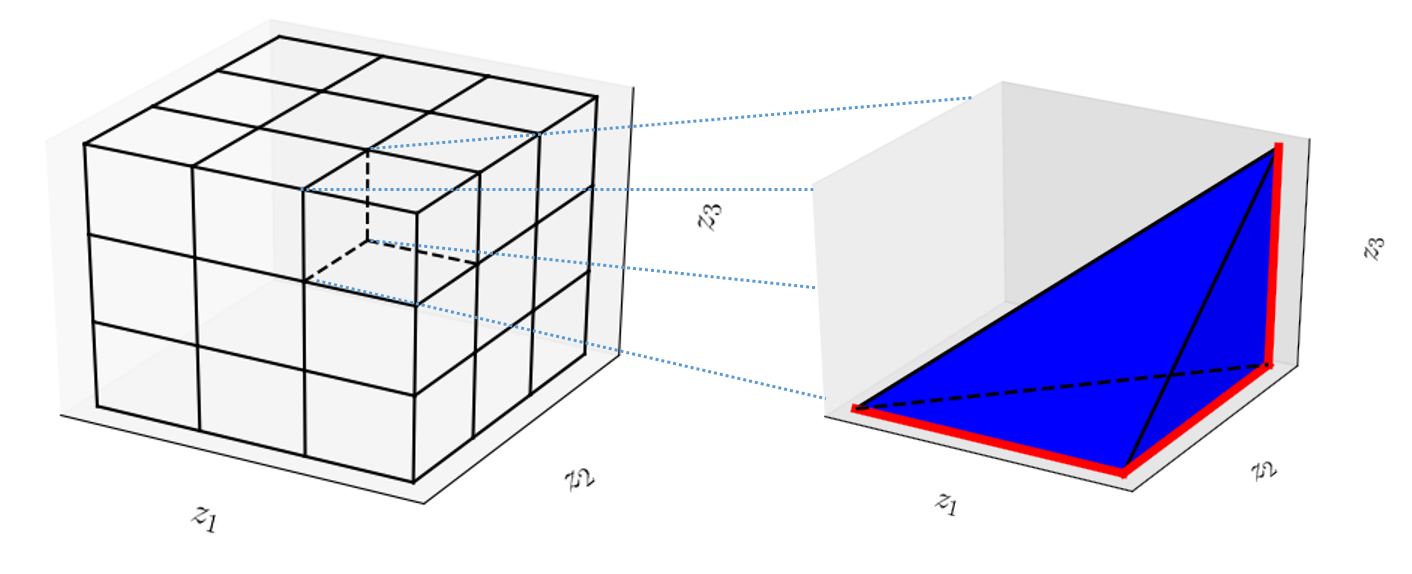}
\caption{Left: The largest hypercube $H$ is defined in $\mathbb{R}^3$ from the minimum and maximum bounds of the three variables $z_1$, $z_2$, and $z_3$. Using four equally spaced breakpoints for each variable, $H$ is further subdivided into $3^3$ smaller cubes. Right: A piecewise 4-simplex inside a subcube is defined by four vertices which constitute a path from the origin of the subcube to the opposite vertex of the subcube as shown by the bold red lines.} \label{Fig:Hypercube2}
\end{figure}

Mathematically, the piecewise linear decomposition of a nonlinear function may be stated as follows, with $f(\mathbf{z})$ as an example and $\hat{f}(\mathbf{z})$ its piecewise linear approximation,
\begin{align}
\hat{f}(\mathbf{z}) = \hat{f}_s(\mathbf{z}), \mathrlap{\enspace \text{if} \enspace \mathbf{z} \in C_s,\, \forall s \in S,\, \forall \mathbf{z} \in H. }
\end{align}
Here, $\mathbf{z}$ is a vector of $d$ variables $\{z_k: z_k^{\text{min}}\leq z_k \leq z_k^{\text{max}}, k \in \{1,\dots,d\}\}$ in $\mathbb{R}^d$. The set $C_s$ is the set of all points in $\mathbb{R}^d$ such that they lie within the simplex of index $s$. The set $S$ is the set of indices of all piecewise $d$-simplices existing in the hyperrectangle $H$ of dimension $d$. The hyperrectangle $H$ is defined for any $\mathbf{z}: z_k \in [z_k^{\text{min}}, z_k^{\text{max}}], \forall k$. In other words, it is the hyperrectangle bounded by the set of vertices $V$ in $\mathbb{R}^d$ such that $V$ contains all $2^d$ unique combinations of the lower and upper bounds of every variable in $\mathbf{z}$. Each piecewise hyperplanar function $\hat{f}_s(\mathbf{z})$ can be expressed as
\begin{align}
\hat{f}_s(\mathbf{z}) = \bm{\nu}_s + \bm{\rho}_s^\mathbf{T}\mathbf{z}. 
\end{align}
Examples of piecewise approximations in one and two dimensions are graphed in Figure \ref{Fig:Piecewise}.
\begin{figure}[h!]
\centerline{
\begin{subfigure}[t]{0.5\textwidth}
\centering
\includegraphics[width=1.0\textwidth]{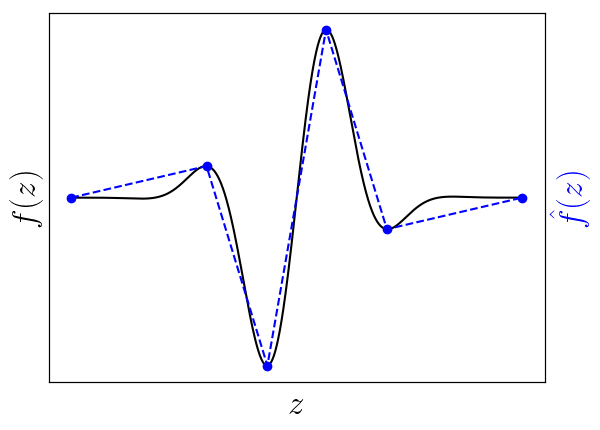}
\caption{}
 \label{Fig:PiecewiseLinear}
\end{subfigure}
~
\begin{subfigure}[t]{0.55\textwidth}
\centering
\includegraphics[width=1.0\textwidth]{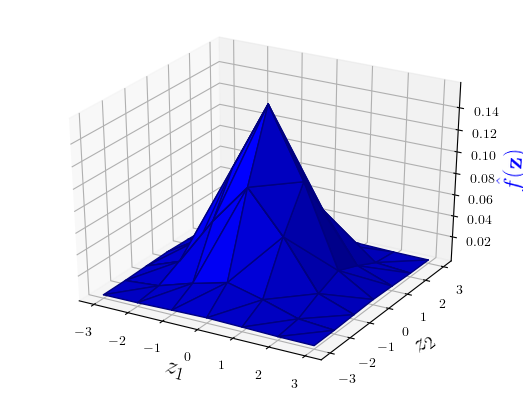}
\caption{}
 \label{Fig:PiecewisePlanar}
\end{subfigure}
}
\caption{(\protect\subref{Fig:PiecewiseLinear}) Piecewise approximation of a nonlinear function $f(z)$ in one dimension. (\protect\subref{Fig:PiecewisePlanar}) Piecewise planar approximation in two dimensions of the multivariate standard normal $f(\mathbf{z})=\exp [-(z_1^2+z_2^2)/2]/2 \pi $}
\label{Fig:Piecewise}
\end{figure}
\label{Fig:Piecewise}
In order to partition the hyperrectangle $H$  into piecewise $d$-simplices, the number and location of breakpoints for each variable $z_k$ in $\mathbf{z}$ must first be decided. Let the breakpoints of $z_k$ be the set of points $\{b_k^l: b_k^0=z_k^{\text{min}}, b_k^{L_k}=z_k^{\text{max}}, l\in \{0,\dots,L_k\}\}$, with $L_k$ the number of piecewise segments and $L_k+1$ the number of breakpoints for the variable $z_k$. The breakpoints of all variables partition $H$ into subspaces occupied by smaller $d$-dimensional hyperrectangles $H_{i}, \forall i \in \{1,\dots,n\}$ such that the interior of any hyperrectangle $H_i$ does not contain a breakpoint as a coordinate. It follows that the number of smaller hyperrectangles or subrectangles is $n=\prod_k L_k$. Each subrectangle $H_i$ can then be further partitioned into piecewise $d$-simplices by using a triangulation scheme. A standard way to triangulate for a $d$-orthotope $H_i$ is for each simplex to be formed from a unique combination of $d$ vertices $V_s$ such that the vertices in $V_s$ constitute a path from the $d$-orthotope that has the origin $(z_1^{\zeta_1}, \dots, z_d^{\zeta_d}): \zeta_k\in \{ 0, \dots, L_k-1 \}$ to the opposing vertex $(z_1^{\zeta_1+1}, \dots, z_d^{\zeta_d+1})$ by trasversing one nearest neighbor vertex at a time. As a result, the number of $d$-simplices for each subrectangle is $m=d!$ meaning that the number of $d$-simplices for $H$ is $mn=d!\prod_k L_k$. An example of a deconstruction of $H$ in three dimensions is shown in Figure \ref{Fig:Hypercube2}. 

\subsection{Modeling piecewise linear approximations for mixed-integer programming}
There are three main schemes for modeling piecewise planar functions in the formulation of MILP problems \cite{Vielma1}. These are the convex combination models, the incremental model, and the multiple choice model. The scheme used in this disclosure belongs to the multiple choice model \cite{Jeroslow1} owing to their relative speed to the other models when the number of simplices is small \cite{Vielma1}. 

In the multiple choice model, we first introduce the binary variables $\{ \mu_{ij}:  i \in \{1,\dots, n\}, j \in \{1,\dots , m\}\}$ with cardinality $mn$. These variables indicate the simplex to pick out of all $mn$ $d$-simplices. Next, a second set of continuous variables $\{ z_k^{ij}\}$ with cardinality $mnd$ is introduced corresponding to $\{ \mu_{ij} \}$ that indicates the value of $z_k$ should a simplex with binary variable $\mu_{ij}$ be chosen. Lastly, we introduce the set $\{\kappa_k^j\}$ where each $\kappa_k^j$ indicates the index of the order of the step taken for variable $x_k^{ij}: \forall i$ from the origin to the opposite vertex for the $j$th simplex. Taken together, the new variables must satisfy the following constraints
\begin{subequations}
\begin{flalign}
&& \sum_{i,j} z_k^{ij} &= z_k, &\forall k\\
&& \sum_{i,j} \mu_{ij} &= 1, & \\
&& b_k^{l} \mu_{ij} \leq z_k^{ij} &\leq b_k^l\mu_{ij} + \frac{b_{k}^{l+1}-b_{k}^{l}}{b_{k'-1}^{l+1}-b_{k'-1}^{l}}(z_{k'-1}^{ij} - b_{k'-1}^{l}\mu_{ij}), & k'=\kappa_k^j, k' \neq 0, \forall k,i,j\\
&& b_k^{l} \mu_{ij} \leq z_k^{ij} &\leq b_k^{l_1}\mu_{ij}, & k'=\kappa_k^j, k' = 0, \forall k, i, j,
\end{flalign}
\end{subequations}
where the breakpoint index $l(i,k)$ for the variable $z_k$ in the $j$th simplex of the $i$th subrectangle corresponds to the lower limit of $z_k$ within the subrectangle. We may then write $\hat{f}(\mathbf{z})$ in terms of the decomposed variables which represent $f(\mathbf{z})$ approximated using affine hyperplanes on the triangulated domain space
\begin{subequations}
\begin{align}
\hat{f}(\mathbf{z}) = \sum_{i,j} \left[ \mu_{ij} f(\mathbf{b}_0^{ij}) + \sum_k \frac{f(\bm{\beta}_k^{ij}) - f(\mathbf{b}_k^{ij})}{b_k^{l+1}-b_k^{l}}(z_k^{ij}-\mu_{ij}b_k^{l}) \right] ,
\end{align}
\begin{flalign}
\text{where} && & & \nonumber \\ 
&& \mathbf{b}_k^{ij}&=(b_1^{l(i,1) + \lambda(i,j,k,1)}, \dots, b_d^{l(i,d) + \lambda(i,j,k,d)}), & \\
&& \lambda(i,j,k,\alpha) &= 1, & \text{if } \kappa_\alpha^j < \kappa_k^j \\
&& \lambda(i,j,k,\alpha) &= 0, & \text{if } \kappa_\alpha^j \geq \kappa_k^j \\
&& \bm{\beta}_k^{ij}&=(b_1^{l(i,1) + \lambda'(i,j,k,1)}, \dots, b_d^{l(i,d) + \lambda'(i,j,k,d)}), & \\
&& \lambda'(i,j,k,\alpha) &= 1, & \text{if } \kappa_\alpha^j \leq \kappa_k^j \\
&& \lambda'(i,j,k,\alpha) &= 0, & \text{if } \kappa_\alpha^j > \kappa_k^j \\
&& \mathbf{b}_0^{ij}&=(b_1^{l(i,1)}, \dots, b_d^{l(i,d)}). & 
\end{flalign}
\end{subequations}

\section{The SPPA algorithm}
The steps of the algorithm are as follows:
\begin{enumerate}
\item An approximation of the nonlinear optimization problem into a linear optimization problem with piecewise linear segments is performed. The number of piecewise segments for each optimization variable involved in a nonlinear expression is defined by the parameter \texttt{initial\_n\_pieces} (i.e. the number of breakpoints is \texttt{initial\_n\_pieces}+1). 
\item \label{Alg:Solution} A solution of the piecewise linear problem is obtained with an MILP solver. 
\item \label{Alg:Contract} With the solution of the piecewise linear optimization problem, the widths of the bounds of optimization variables involved in nonlinear expressions are contracted proportionally by a factor of \texttt{contract\_frac}$<1$. The newly contracted bounds of an optimization variable are centered about its solution. If the newly contracted bounds of an optimization variable violate the old bounds such that the new lower (upper) bound is smaller than the old lower (upper) bound, the newly contracted bounds are translated positively (negatively) such that the lower (upper) bound is at the position of the lower (upper) bound of the previous bounds. 
\item \label{Alg:Adjust} The piecewise linear optimization problem is reconstructed using the contracted bounds where each optimization variable involved in one or more nonlinear expression has \texttt{n\_pieces} number of piecewise segments. 
\item Step \ref{Alg:Solution} can be repeated again with the reconstructed optimization problem and the bounds can be contracted further with Step \ref{Alg:Contract} followed by another reconstruction of the piecewise linear optimization problem. In this way, Steps \ref{Alg:Solution}, \ref{Alg:Contract}, and \ref{Alg:Adjust} are repeated until one or more desired termination criteria have been reached. 
\end{enumerate}

\section{Results}
In this section, we will present the optimization results for the two-dimensional Rosenbrock function (NLP), the two-dimensional Rastrigin function (NLP), the two-dimensional Ackley function (NLP), the Eggholder function (NLP), a spring design problem (MINLP) \cite{Sandgren1,Leyffer1}, and a cyclic scheduling problem (MINLP) \cite{Jain1,Leyffer3}. 

\renewcommand{\arraystretch}{1.3}
\begin{table}[h!] 
\begin{tabular}{l l l l l l l}
\toprule
                                                                       & \multicolumn{2}{c}{Objective} & \phantom{as} & \multicolumn{2}{c}{Segments} & \\
\cmidrule{2-3} \cmidrule{5-6}
Problem & Found & Optimal & & initial\_n\_pieces & n\_pieces & CPU time \\
\midrule
Rosenbrock               &                     6.13E-6                      &           0                            &  &     4        &    4    &   1.5s   \\
Rastrigin                      &                             0                              &                0                        & &     6            &  3      &    0.5s  \\
Ackley                            &                      2.7E-6                              &                     0                  &  &       3        &     3   &    1.8s   \\
Eggholder                     &                    -959.6407                              &                    -959.6407             &  &    35        &      3   &    4.4s \\
Spring design                 &                   0.84625                          &                     0.84625                    &  &     3        &     3  &   95s   \\
\midrule
Cyclic scheduling \\
K=4, \, $\epsilon$=0                             &                 -165920.1             &     -165920.1        &  & 4  & 3    &   18 mins    \\
K=4, \, $\epsilon$=0.01                          &                  -165398.7            &        -165398.7            &   &    4    &   3   & 16 mins  \\
K=10, $\epsilon$=0                              &                    -166322.0           &        -166322.0      &    &    4   &     3  &  16 mins  \\
K=10, $\epsilon$=0.01                      &                  -166102.0                  &      -166102.0        &  &    4     &    3 &    13 mins   \\
\bottomrule
\end{tabular}
\caption{Summary of results obtained and parameters used for the test problems. }
\label{Table:Summary}
\end{table}
For the cyclic scheduling problem, $K$ is a parameter of the optimization problem while $\epsilon$ is a small constant used in the objective function to avoid a divide-by-zero situation. Because of the derivative-free nature of piecewise linear approximations, $\epsilon$ can be set to zero so we also present an optimal solution found by SPPA that has not been reported in the literature. We compare it against the solution obtained by BONMIN from a slight MINLP reformulation of the problem with extra binary variables in order to avoid having $\epsilon$. Also, modeling the problem directly with SPPA creates MILP problems which are intractable. This problem may be circumvented by solving simpler subproblems where the optimization variable $T_{\text{cycle}}$ is kept constant. In this way, the optimal solution may be obtained by performing a binary search for $T_{\text{cycle}}$ over the smaller subproblems since the optimization problem is pseudo-convex. 

A summary of the results obtained is provided in Table  \ref{Table:Summary}. A laptop with an Intel i7-6820HQ CPU was used for computation while CPLEX was used for solving the MILP problems. 

\section{Discussion}
The intuition behind SPPA is that for most optimization problems, there isn't a need to obtain very fine approximations for regions of low interest where a coarser approximation would very likely be indicative of the feasibility and optimality of a domain subset. Furthermore, for an optimization landscape that exhibits a certain degree of self-similarity in a region of high interest, an iteratively and exponentially shrinking piecewise approximation with a small number of integer variables would be much more prefarable than using more breakpoints and integer variables for obtaining a solution of higher precision. This is because the bounded domain shrinks exponentially with the number of iterations in SPPA whereas more integer variables could easily lead to a problem scaling exponentially in difficulty. 

Feasibility with SPPA can only be guaranteed with linear constraints. Global and even local optimality for SPPA cannot be guaranteed although our results show that the algorithm works quite well in practice. Our results also show that a small number of segments and integer variables can be sufficient in order for SPPA to converge on the optimal solution which is important for limiting the additional computational complexity that comes with piecewise linear approximations. 

Altogether, SPPA extends and improves on piecewise linear programming methods by offering a simple and iterative scheme by which to converge on optimal solutions using only a small number of breakpoints and integer variables. 

\paragraph{Software} The code for the SPPA solver is available as a Python package called \texttt{sppa} and is available at \url{https://pypi.org/project/sppa/} and \url{https://github.com/jamespltan/sppa}. 

\paragraph{Funding} This research did not receive any specific grant from funding agencies in the public, commercial, or not-for-profit sectors. 

\paragraph{Declaration of Interest} The author declares no competing interests, financial or otherwise. 

\bibliographystyle{elsarticle-num-names}
\bibliography{References}

\end{document}